\documentclass[10pt,a4paper]{article}
\usepackage{amsmath}
\usepackage{amssymb}
\usepackage{amsthm}
\usepackage{amsfonts}
\usepackage{amscd}
\usepackage[mathscr]{eucal}
\usepackage{multirow}
\newcommand{\Z} {{\mathbb  Z}}
\newcommand{\Q}{{\mathbb  Q}}
\newcommand{\F}{{\mathbb  F}}

\newcommand{\R} {{\mathbb R}}

\begin{document}
\parindent  25pt
\textwidth  15cm    \textheight  23cm \evensidemargin -0.06cm
\oddsidemargin -0.01cm
\linespread{1.5} 

\title{ On $ p-$Rationality of Cubic and Quartic Number Fields }
\author{\mbox{} {Hang Li$^1$ and
Derong Qiu$^2$}\thanks{ Corresponding author.
E-mail addresses: 12331008@mail.sustech.edu.cn (H.Li); \
derong@mail.cnu.edu.cn (D.Qiu) } \\
(1.Department of Mathematics, Southern University of Science and Technology,\\
Shenzhen 518055,P.R.China\\
2.School of Mathematical Sciences, Capital Normal University, \\
 Beijing 100048, P.R.China) }

\date{}
\maketitle
\parindent  24pt
\parskip  0pt

\par   \vskip 0.6cm

{\bf Abstract} \quad In this paper, a new criterion is given to determine
the $p-$rationality of some complex cubic number fields in terms of $ p-$divisibility
of certain terms of a third-order recurrence sequence, several illustrated examples are constructed,
the relations between generalized $ abc-$conjecture and the $p-$rationality are discussed,
from which some explicit fields satisfying Greenberg's Generalized Conjecture (GGC, for short)
are obtained.

\par  \vskip  0.6 cm

\par  \vskip  0.3 cm

\par     \vskip  1 cm

\hspace{-0.6cm}{\bf 1. \ Introduction and Main Results}

      \par \vskip 0.6 cm

Let $ F $ be a number field, $ p $ a prime number. Let $ S_{p} $ denote the set of
all places of $ F $ above $ p, \ M_{S_{p}}(F) $ the maximal pro-$p$ extension unramified
outside $ S_{p}, $ and $ F_{S_{p}} $ the maximal Abelian pro-$p$ extension unramified
outside $ S_{p}. $ Write $ G_{S_{p}}(F): = \text{Gal}(M_{S_{p}}(F)/F), $ then
$ G_{S_{p}}(F)^{ab} =  \text{Gal}(F_{S_{p}}/F). $ In the following, we denote by
$ \mathcal{T} _{F} = \text{tor}(G_{S_{p}}(F)^{ab}) $ the maximal $ \Z_{p}-$torsion
submodule of $ G_{S_{p}}(F)^{ab}. $ \\
$ F $ is called $p-$rational if $ G_{S_{p}}(F) $ is a free
pro-$p$ group (see [Mov1, 2]).  \\
It is known ([Gr2],section 3) that a number field $ F $ is $p-$rational if and only if the
following two conditions hold \\
(1) $  \text{rank}_{\Z_{p}}(G_{S_{p}}(F)^{ab}) = r_{2}(F) + 1; $ \\
(2) \ $ G_{S_{p}}(F)^{ab} $ is a torsion free $ \Z_{p}-$module, i.e.,
$ \mathcal{T} _{F} = \{1 \}. $  \\
Condition (1) above is equivalent to say that Leopoldt's conjecture holds for $ (F, p), $
while condition (2) means that $ F_{S_{p}} = \widetilde{F}, $ the composite field of
all the $ \Z_{p}-$extension of $ F. $ \\
The $ p-$rationality of number fields has been extensively studied, one of its important
application is on constructing non-abelian number fields satisfying Leopoldt's conjecture
(e.g. [Mov1, 2]). Besides, much questions about it are still open, such as \\
{\bf Conjecture 1.1} (Gras [G2, Conjecture 8.11]). {\it Every number field $ F $ is
$ p-$rational if $ p $ is sufficiently large. }  \\
{\bf Conjecture 1.2} (Greenberg [Gr2]). {\it For any odd prime number $ p $ and positive integer
$ t, $ there exists a Galois extension $ F $ of the rational number field $ \Q $
such that $ F $ is $ p-$rational and $ \text{Gal}(F/ \Q) \cong (\Z / 2 \Z)^{t}. $ } \\
When $ F $ is a quadratic number field or some of their composite, much progress have
been made on the Greenberg's conjecture above (see e.g., [AB],[BM], [Ko], [CLS]). \\
There are a series of work on determining the $ p-$rationality of totally real number fields
and quadratic imaginary number fields (see e.g., [Gr2], [Li]). In this paper, we will study
similar questions on some cases that $ F $ is neither totally real nor quadratic, more
precisely, we consider the $ p-$rationality of complex cubic and pure imaginary quartic Number
Fields. Our main results are the following
\par \vskip 0.2 cm

{\bf Theorem 1.3.} \ {\it Let $ K \subset \R $ be a complex cubic  number field,
$ \varepsilon > 0 $ be its fundamental unit, $ f(x) = x^{3} -a_{2} x^{2} - a_{1}x - a_{0}
\in \Z[x] $ be the minimal polynomial of $ \varepsilon $ over $ \Q. $
Let $ p $ be an odd prime number with $ p \nmid d(f) $ (the discriminant of $ f(x) $).
Define the arithmetic third order sequence of rational integers \\
$ F_{n+3} = a_{2} F_{n+2} + a_{1} F_{n+1} + a_{0} F_{n} \ (n \geq 0),  \
F_{0} = F_{1} = 0, F_{2} = 1. $  \\
Then we have the following conclusion:  \\
(1) \ if $ p $ splits completely over $ K / \Q, \ F_{p-1} \not\equiv 0 (\text{mod} \ p^{2}), $
then there is a prime ideal $ \mathfrak{p} $ of $ K $ lying over $ p $ such that
$ \varepsilon^{p - 1} \not\equiv 1 (\text{mod} \ \mathfrak{p}^{2}). $  \\
(2) \ if $ p \mathcal{O}_{K} = \mathfrak{p}_{1}\mathfrak{p}_{2} $ for two different
prime ideals in $ K $ with $ f(\mathfrak{p}_{1}/p) = 1 $ and $ f(\mathfrak{p}_{2}/p) = 2, $
then \\
$  F_{p^{2}-1} \not\equiv 0 (\text{mod} \ p^{2}) \Rightarrow
\varepsilon^{p - 1} \not\equiv 1 (\text{mod} \ \mathfrak{p}_{1}^{2}) $ or
$ \varepsilon^{p^{2} - 1} \not\equiv 1 (\text{mod} \ \mathfrak{p}_{2}^{2}). $ \\
(3) \ if $ p $ is inert over $ K / \Q, $ i.e., $ p \mathcal{O}_{K} = \mathfrak{p} $
is a prime ideal in $ K, $ then  \\
$  F_{p^{3}-1} \not\equiv 0 (\text{mod} \ p^{2}) \Rightarrow
\varepsilon^{p^{3} - 1} \not\equiv 1 (\text{mod} \ \mathfrak{p}^{2}). $ \\
Furthermore,if $p$ satisifes one of (1)-(3) and $p\nmid h(K),$ where $h(K)$ is the class number
of $K,$ then $K$ is $p-rational.$ }
\par \vskip 0.2 cm

{\bf Theorem 1.4.} \ {\it Let $ K $ be a complex cubic or pure imaginary quartic number field.
If the generalized $ abc-$conjecture holds for $ K, $ then as $ x \rightarrow \infty , $ we have
$$ \sharp \{\text{prime numbers} \ p \leq x : \ K \ \text{is} \ p-\text{rational} \} \geq c \log x $$
for some constant $ c > 0 $ depending only on $ K. $ }
\par \vskip 0.2 cm

In the following, for a number field $ F, $ as usual, we denote by $ r_{1} $ and
$ r_{2} $ the number of real and complex infinite places, respectively.
So $ r_{1} + 2 r_{2} = [F : \Q]. $ Let $ \mathcal{O}_{F}$ be the ring of integers of $ F $
and $ E_{F} = \mathcal{O}_{F}^{\times } $ denote the group of units. For any place $ v $
of $ F, $ let $ F_{v} $ denote the completion of $ F $ at $ v, $ and $ U_{v}^{(1)} $
denote the group of principal units of $ F_{v}. $ Also, we let $ h(F) $ denote the class
number of $ F. $ \\
In general, for a finitely generated $ \Z_{p}-$module $ M, $ let $ \text{tor}(M) $
denote its maximal torsion $ \Z_{p}-$submodule.

\par     \vskip  0.6 cm

\hspace{-0.6cm}{\bf 2. \ A criteria of $ p-$rationality}

      \par \vskip 0.6 cm

Let $ K $ be a complex cubic number field, then $ [K : \Q] = 3 $ and $ r_{1} = r_{2} = 1. $
Let $ \varepsilon > 1 $ be its fundamental unit. For any prime number $ p, $ it is easy to
see that the Leopoldt's conjecture ([Ka]) holds for $ (K, p). $
\par \vskip 0.2 cm

{\bf Proposition 2.1.} ([NG, Prop.6.1],[BM,Prop.2.3]) {\it A number field $F$ is $p-$rational
if and only if the following three conditions are satisfied:\\
(1) The $p-$ Hilbert class field $H_{p}(F)\subseteq \widetilde{F};$\\
(2) The natural map $\mu_{p}(F)\rightarrow \oplus_{v|p}\mu_{p}(F_{v})$ is an isomorphism and \\
(3)The natural map $U_{F}/p\rightarrow \oplus_{v|p}U_{v}/p$ is injective.
} \\
Here $\mu_{p}(F)$ consists of the $p-$th roots of unity of any field $ F, \ U_{F}$
is the group of global units and $U_{v}$ the group of local units in $F_{v}.$

{\bf  Proposition 2.2.} \ {\it Let $ K $ be a complex cubic or pure imaginary quartic number field,
$ p $ a prime number (when $ p = 3, $ assume it is unramified in $ K / \Q, $  when
$ p = 5 $ and $ K $ is quartic, assume $ 5 $ is not totally ramifed in $ K $).  \\
Then $ K $ is $ p-$rational if and only if the following two conditions hold: \\
(1) \  $ H_{p}(K) \subseteq \widetilde{K}, $ where $ H_{p}(K) $ is the Hilbert $ p-$class
field of $ K, $ and $ \widetilde{K} $ is the composite of
all $ \Z_{p}-$extension of $ K. $  \\
(2) \ There is a prime ideal $ \mathfrak{p} $ of $ K $ lying over $ p $ such that
$ \varepsilon^{p^{f} - 1} \not\equiv 1 (\text{mod} \ \mathfrak{p}^{e+1}) $  for
any fundamental unit $ \varepsilon $ of $ K, $  where $ e $ and $ f $
are the ramification index and residue degree of $ \mathfrak{p} $ in $ K /\Q, $ respectively. }
\par \vskip 0.2 cm

{\bf  Proof.} \ We only prove the case that $ K $ is a complex
cubic number field, the other case is similar.
 By above Prop.2.1, we only need check that the number field
$K$ satisfies conditions (2) and (3) in the above Prop.2.1. If $p>5$ is an odd prime or
$p=3$ is unramified in $K/\Q$, we always have $\mu_{p}(K)=1,\oplus_{v|p} \mu_{p}(K_{v})=1.$
Hence the condition (2) holds. Now we claim that the condition (3) of above Prop.2.1 is
equivalent to condition (2) of this proposition. In fact, the factor group $U_{K}/U_{K}^{p}$ is
generated by the class of the fundamental unit $\varepsilon$. We denote by $U_{v}^{(i)}$ the higher
unit groups of the local field $K_{v}$. Then we have $U_{v}/U_{v}^{p}=U_{v}^{(1)}/(U_{v}^{(1)})^p$
since $\mu_{p}(K_{v})=\{1\}.$ Let $\mathfrak{p}_{v}$ be the prime correspond with $v$, we have
$\frac{e_{v}}{p-1}<1,$ where $e_{v}$ is the ramification index of $\mathfrak{p}$ in $K/\Q$. By
Proposition 9,Chapter 14 of [S],the map $x\mapsto x^p$ is an isomorphism of $U_{v}^{(1)}$ onto
$U_{v}^{(1+e_{v})}$.Hence the natural map $U_{K}/U_{K}^p\rightarrow \oplus_{v|p}U_{v}/U_{v}^{p}$
is injective if and only if there is a prime ideal $\mathfrak{p}$ of $K$ lying over
 $p$ such that $ 1 \neq \varepsilon \in U_{\mathfrak{p}}^{(1)}/(U_{\mathfrak{p}}^{(1)})^{p} =
 U_{\mathfrak{p}}^{(1)}/U_{\mathfrak{p}}^{(1+e_\mathfrak{p})}$. It is easy to see that
 $ 1 \neq \varepsilon \in U_{\mathfrak{p}}^{(1)}/U_{\mathfrak{p}}^{(1+e_\mathfrak{p})}$ is
 equivalent to the global condition $  \varepsilon^{p^{f} - 1} \not\equiv 1 (\text{mod} \ \mathfrak{p}^{e+1}), $
 where $ e $ and $ f $  are the ramification index and residue degree of $ \mathfrak{p} $ in $ K /\Q, $
 respectively. The proof of this proposition is completed.
\quad $ \Box $

\par     \vskip  0.6 cm

\hspace{-0.6cm}{\bf 3. \ Proof of Theorem 1.3}

      \par \vskip 0.6 cm

{\bf Proof of Theorem 1.3.} \ We have $ K = \Q(\varepsilon ), d_{K}\mid d(f), $ where
$ d_{K} $ is the discriminant of $ K. $ So by assumption $ p $ is unramified in $ K. $
Let $ \sigma _{1} = \text{id}, \sigma _{2}, \sigma _{3} $ be the three imbedding
of $ K $ into $ \overline{\Q} $ over $ \Q. $ Then $ f(x) = (x - \varepsilon )
(x - \sigma _{2} ( \varepsilon ))(x - \sigma _{3}(\varepsilon )).  $ By the general
term formula,
$$ \begin{array}{l} F_{n} = A \varepsilon ^{n} + B \sigma _{2} ( \varepsilon )^{n} +
C \sigma _{3} ( \varepsilon )^{n}, \ \text{where}  \\
A = [(\sigma _{1} ( \varepsilon ) - \sigma _{2} ( \varepsilon ))
(\sigma _{1} ( \varepsilon ) - \sigma _{3} ( \varepsilon ))]^{-1}, \\
B = [(\sigma _{2} ( \varepsilon ) - \sigma _{1} ( \varepsilon ))
(\sigma _{2} ( \varepsilon ) - \sigma _{3} ( \varepsilon ))]^{-1}, \\
C = [(\sigma _{3} ( \varepsilon ) - \sigma _{1} ( \varepsilon ))
(\sigma _{3} ( \varepsilon ) - \sigma _{2} ( \varepsilon ))]^{-1}. \end{array} $$
Notice that $ A + B + C = 0, $ so $ F_{n} = B (\sigma _{2} ( \varepsilon )^{n} -
\varepsilon ^{n}) + C (\sigma _{3} ( \varepsilon )^{n} - \varepsilon ^{n}), $
and so
$$ \begin{array}{l} (\sigma _{2} ( \varepsilon ) - \sigma _{1} ( \varepsilon ))
(\sigma _{2} ( \varepsilon ) - \sigma _{3} ( \varepsilon ))
(\sigma _{3} ( \varepsilon ) - \sigma _{1} ( \varepsilon )) F_{n} \\
= (\sigma _{3} ( \varepsilon ) - \sigma _{1} ( \varepsilon ))
(\sigma _{2} ( \varepsilon )^{n} - \varepsilon ^{n}) -
(\sigma _{2} ( \varepsilon ) - \sigma _{1} ( \varepsilon ))
(\sigma _{3} ( \varepsilon )^{n} - \varepsilon ^{n}).  \end{array} $$
Let $ D = (\sigma _{2} ( \varepsilon ) - \sigma _{1} ( \varepsilon ))
(\sigma _{2} ( \varepsilon ) - \sigma _{3} ( \varepsilon ))
(\sigma _{3} ( \varepsilon ) - \sigma _{1} ( \varepsilon )), $  then
$ D^{2} = d(f). $ Now we come to show the results (1), (2) and (3) of
the Theorem. Our way of proof is via contradiction. \\
(1) \ Let $ L $ be the normal closure of $ K/\Q, $ each $ \sigma _{i}
(i = 1, 2, 3) $ can be extended on $ L, $ and denoted by the same notations. If we
suppose that $ \varepsilon^{p - 1} \equiv 1 (\text{mod} \ p^{2}\mathcal{O}_{K}), $
then $ \varepsilon^{p - 1} \equiv 1 (\text{mod} \ p^{2}\mathcal{O}_{L}). $
So $ \sigma _{i}( \varepsilon^{p - 1}) \equiv 1 (\text{mod} \
\sigma _{i}(p^{2}\mathcal{O}_{L}) = p^{2}\mathcal{O}_{L})
(i = 1, 2, 3). $ Hence $ \varepsilon^{p - 1} \equiv \sigma _{2}( \varepsilon^{p - 1})
\equiv \sigma _{3}( \varepsilon^{p - 1}) (\text{mod} \ p^{2}\mathcal{O}_{L}), $
which implies that $ D F_{p-1} \in p^{2}\mathcal{O}_{L}, $ and so
$ F_{p-1} \in p^{2}\mathcal{O}_{L} $ because $ p \nmid d(f) = D^{2}. $
Therefore $ p^{2} \mid F_{p-1} $ because $ F_{p-1} \in \Z $ is an integer. This
contradicts to the assumption. This proves (1). \\
(2) and (3) can be similarly done. If $p$ further satisfies $p\nmid h(K),$
both (1) and (2) of Proposition 2.2 are satisfied. so $K$ is $p-$rational.    \quad $ \Box $
\par     \vskip  0.6 cm

\hspace{-0.6cm}{\bf 4. \ Proof of Theorem 1.4}
      \par \vskip 0.6 cm

{\bf Proof of Theorem 1.4.} \ Our proof is in the same spirit of ([MR], second part).
For the generalized $ abc-$conjecture, see e.g., [Ic], [MR]. Let $ \varepsilon $ be
the fundamental unit of $ K. $ Then $ \exists \ n_{0} \in \Z_{> 0} $
such that for every integer $ n \geq n_{0} $ with $ \varphi (n) \geq n / 2 $ ($ \varphi $
is the Euler function), there is a prime ideal $ \mathfrak{p}_{n} $ of $ \mathcal{O}_{K} $
satisfying \\
(1) \ $ \varepsilon ^{N\mathfrak{p}_{n} - 1} \not\equiv 1 (\text{mod} \ \mathfrak{p}_{n}^{2}); $ \\
(2) \ $ \mathfrak{p}_{n} \neq \mathfrak{p}_{n^{\prime }} $ if $ n \neq n^{\prime }\geq n_{0}; $ \\
(3) \ there is a real number $ r > 0 $ depending only on $ K, $ such that $ N \mathfrak{p}_{n}
\leq r^{n}. $  \\
Here $ N \mathfrak{p}_{n} = \sharp (\mathcal{O}_{K} / \mathfrak{p}_{n}) $ is the norm of ideal. \\
In fact,  this can be proved by using the same method as Proposition 2.3 of [MR]. We only outline the main steps,
for details, see Section 2 of [MR]. We know that $\varepsilon$ has no conjugate on the unit circle. So by Lemma 2.1 of [MR],
there exists some $ k \in \Z_{> 0}$ such that $$ |N(\Phi_{n}(\varepsilon^{k}))| \geq \exp (cn) $$
for $n$ such that $\varphi(n)\geq \frac{1}{2}n, $ where $ c>0 $ is a constant depending on $ \varepsilon^{k}, \ N() $ is the
norm map form $K$ to $ \Q $ and $ \Phi_{n}(x) $ is $ n-$th cyclotomic polynomial.
We denote $ u=\varepsilon^{k}. $ Let $ (u^{n}-1) = I_{n}J_{n},(\Phi_{n}(u)) = A_{n}B_{n} $ with $ I_{n} $ and $ J_{n} $
relatively prime and where if $ \mathfrak{p}|I_{n}, $ then $ \mathfrak{p}^{2}\nmid I_{n}, $ and if $ \mathfrak{p}|J_{n}, $
then $ \mathfrak{p}^{2}\mid J_{n}, $ and $A_{n}, B_{n}$ are similar. Applying the generalized $abc$ conjecture to the
triple $ (u^{n}-1,1,u^{n}), $ we obtain  for all $ \varepsilon>0, N(J_{n})\ll_{K,\varepsilon}|N(u^{n} - 1)|^{2\varepsilon}. $
From the above two inequalities and choose $ \beta>1 $ such $ |\sigma_{i}(u)|\leq \beta, $ where
$ \sigma_{i}(i=1,2,\cdots, [K:\Q])$ is an embedding $ \sigma_i: K \hookrightarrow \mathbb{C}, $ then we have $$
N(A_{n})\gg_{K,\varepsilon} \text{exp}(n(c-2m\varepsilon \log \beta)) $$
for every $n$ such that $\varphi(n)\geq \frac{1}{2}n, $ where $ m=[K,\Q]. $
From this, we can get that there exists $ \ n_{0} \in \Z_{> 0} $
such that for every integer $ n \geq n_{0} $ with $ \varphi (n) \geq n / 2, $ there is a prime ideal
$ \mathfrak{p}_{n} $ of $ \mathcal{O}_{K} $ satisfying \\
(1)$^{\prime}$ \ $ \mathfrak{p}_{n} |\Phi_{n}(\varepsilon^{k}) $ and $ \varepsilon ^{kn} \not\equiv 1 (\text{mod} \ \mathfrak{p}_{n}^{2}); $ \\
(2)$^{\prime}$ \ $ \varepsilon^{k} $ is of order $n$ in $ (\mathcal{O}_{k}/\mathfrak{p}_{n})^{\times}; $ \\
(3)$^{\prime}$ \ there is a real number $ r > 0 $ depending only on $ K, $ such that $ N \mathfrak{p}_{n} \leq r^{n}. $ \\
Combining this with Lemma 4 of [Ic], we obtain such $ \mathfrak{p}_{n} $ satisfying
$ \varepsilon ^{N \mathfrak{p}_{n}-1} \not\equiv 1 (\text{mod} \ \mathfrak{p}_{n}^{2}) $
and from (2)$^{\prime}, $ we know $ \mathfrak{p}_{n}\neq \mathfrak{p}_{n'} $ if $ n\neq n'\geq n_{0}. $ \\
Let $ x \geq 1, $ for the above real number $ r, $ take $ n_{1} = [\log_{r} x]. $ Let $ x $
be sufficiently large such that $ n_{1} \geq n_{0}. $ Then for any integer $ n: n_{0} \leq n
\leq n_{1} $ with $ \varphi (n) \geq n /2, $ there exists a prime ideal $ \mathfrak{p}_{n} $
of $ \mathcal{O}_{K} $ such that $ \varepsilon ^{N\mathfrak{p}_{n} - 1} \not\equiv 1 (\text{mod} \
\mathfrak{p}_{n}^{2}). $ Let $ \mathfrak{p}_{n} \cap \Z = p_{n} \Z $ for some prime number $ p_{n}. $
Then $  p_{n} \leq N\mathfrak{p}_{n} \leq r^{n} \leq r^{n_{1}} \leq x. $ So
$ \frac{1}{m} \cdot \sharp \{n: n_{0} \leq n \leq n_{1}, \varphi (n) \geq n / 2 \} \leq
\sharp \{\text{primes} \ p_{n} \leq x: \exists \ \text{prime ideal} \ \mathfrak{p}_{n} \lhd
\mathcal{O}_{K}, \mathfrak{p}_{n} \mid p_{n},
\varepsilon ^{N\mathfrak{p}_{n} - 1} \not\equiv 1 (\text{mod} \ \mathfrak{p}_{n}^{2}) \}. $
From Silverman's lower bound (see Lemma 6 of [Sil]) $$
   \sharp \{n\leq x:\varphi(n)\geq \frac{n}{2} \}\geq (\frac{6}{\pi^2}-\frac{1}{2})x +O(\log x)
    $$ and $ n_{1}=[\text{log}_{r}x], $
    we get \[\begin{split}
       & \sharp \{n:n_{0}\leq n\leq n_{1},\varphi(n)\geq \frac{n}{2} \} \\
       =&\sharp \{n: n\leq  [\log_{r} x],\varphi(n)\geq \frac{n}{2} \}-\sharp\{0\leq m< n_{0}:\varphi(m)\geq \frac{m}{2}\}\\
       =&\sharp \{n: n\leq  \text{log}_{r}x,\varphi(n)\geq \frac{n}{2} \}-\sharp\{0\leq m< n_{0}:\varphi(m)\geq \frac{m}{2}\}\\
       \geq&(\frac{6}{\pi^2}-\frac{1}{2})\log_{r} x +O(\log\log_{r} x)- \sharp\{0\leq m< n_{0}:\varphi(m)\geq \frac{m}{2}\}\\
        =&\frac{1}{\log r}(\frac{6}{\pi^2}-\frac{1}{2})\log x +O(\log \log_{r} x)- \sharp\{0\leq m<n_{0}:\varphi(m)\geq \frac{m}{2}\},\\
    \end{split}
    \]
so as $x\rightarrow \infty,\sharp \{n:n_{0}\leq n\leq n_{1},\varphi(n)\geq \frac{n}{2} \}\geq c\cdot \log x $
for some real number $ c > 0 $ which depends only on $ r, $ hence depends only on $ K.$
Therefore, there exists a real number $ c > 0 $ depending only on $ K, $
such that, as $ x \rightarrow \infty , \ \sharp \{\text{primes} \ p_{n} \leq x: \exists \
\text{prime ideal} \ \mathfrak{p}_{n} \lhd \mathcal{O}_{K}, \mathfrak{p}_{n} \mid p_{n},
\varepsilon ^{N\mathfrak{p}_{n} - 1} \not\equiv 1 (\text{mod} \ \mathfrak{p}_{n}^{2}) \}
\geq c \cdot \log x. $ Let $ A = \{\text{prime} \ p: p \mid h(K) \ \text{or} \
p \ \text{is ramified in} \ K \} $  ($ h(K) $ is the class number of $ K $). Then $\sharp A < \infty . $
The complement set of $ A $ in the set of all primes is
$ A^{c}=\{\text{prime} \ p: p \nmid h(K) \ \text{and} \ p \ \text{is unramified in} \ K  \}. $
By Proposition 2.2 above, we know
$$\{\text{primes} \ p_{n} \leq x: \exists \
\text{prime ideal} \ \mathfrak{p}_{n} \lhd \mathcal{O}_{K}, \mathfrak{p}_{n} \mid p_{n},
\varepsilon ^{N\mathfrak{p}_{n} - 1} \not\equiv 1 (\text{mod} \ \mathfrak{p}_{n}^{2}) \}
\cap A^{c} $$
$$\subseteq \{\text{primes} \ p \leq x :
K \ \text{is} \ p-\text{rational} \}.
$$
So, as $ x \rightarrow \infty , \ \sharp \{\text{primes} \ p \leq x :
K \ \text{is} \ p-\text{rational} \} \geq \sharp (\{\text{primes} \ p_{n} \leq x: \exists \
\text{prime ideal} \ \mathfrak{p}_{n} \lhd \mathcal{O}_{K}, \mathfrak{p}_{n} \mid p_{n},
\varepsilon ^{N\mathfrak{p}_{n} - 1} \not\equiv 1 (\text{mod} \ \mathfrak{p}_{n}^{2}) \}
\setminus A) \geq \sharp \{\text{primes} \ p_{n} \leq x: \exists \
\text{prime ideal} \ \mathfrak{p}_{n} \lhd \mathcal{O}_{K}, \mathfrak{p}_{n} \mid p_{n},
\varepsilon ^{N\mathfrak{p}_{n} - 1} \not\equiv 1 (\text{mod} \ \mathfrak{p}_{n}^{2}) \}
- \sharp A \geq c \log x - \sharp A \geq c^{\prime } \log x . $
The proof of Theorem 1.4 is completed. \quad $ \Box $
\par \vskip 0.2 cm

{\bf Remark.} \ As pointed by the referee, once non-divisibility of the class number $ h(K) $
by $ p $ is known, verifying $ p-$rationality from the recurrence sequence can be quick and
useful. On the other hand, in general, the condition $ p | h(K) $ does not imply $ K $ is
not $ p-$rational.

\par     \vskip  0.6 cm

\hspace{-0.6cm}{\bf 5.\ Examples }

      \par \vskip 0.6 cm

{\bf Example 5.1.} \ Let $ p\geq 5 $ be a prime number, $ K = \Q(\alpha) $ with
$ \alpha = \sqrt[3]{p^{3}-1}. $ Then $ p - \alpha $ is a fundamental unit of
$ K $ (see [St]). Let $ \varepsilon = \frac{1}{p - \alpha} = p^{2} + pa + \alpha ^{2}, $
then $ \varepsilon > 1 $ is also a fundamental unit of $ K. $ The minimal polynomial
of $ \alpha $ on $ \Q $ is $ f(x) = x^{3} + 1 - p^{3} $ with discriminant $ d(f) =
-27 (1 - p^{3})^{2}. $ It can be easily verified that
(see e.g. [N],Chapt 1, Prop.8.3), \\
(R1) \ if $ p \equiv 1 (\text{mod} \ 3), $ then $ p $ splits completely in $ K / \Q, $
i.e., $ p \mathcal{O}_{K} = \mathfrak{p}_{1}\mathfrak{p}_{2}\mathfrak{p}_{3} $ for
three deferent prime ideals in $ K; $ \\
(R2) \ if $ p \equiv 2 (\text{mod} \ 3), $ then $ p \mathcal{O}_{K} =
\mathfrak{p}\mathfrak{q} $ for two different prime ideals in $ K $
with residue degree $ f(\mathfrak{p}/p) = 1,  f(\mathfrak{q}/p) = 2. $ \\
For (R1), $ f = f(\mathfrak{p}_{i}/p) = 1 (i = 1, 2, 3). $ We have
$ \varepsilon ^{p -1} = (p^{2} + pa + \alpha ^{2})^{p-1} = (p(p+\alpha) + \alpha ^{2})^{p-1}
\equiv \alpha ^{2p-2} + (p-1)\alpha ^{2p-4}p(p + \alpha ) \equiv
(\alpha ^{3})^{\frac{2(p-1)}{3}} + (p-1)\alpha ^{2p-4}p \alpha
\equiv (p^{3} - 1)^{\frac{2(p-1)}{3}} - p \alpha ^{2p-3} \equiv
1 - p \alpha ^{2p-3} \not\equiv 1 (\text{mod} \ p^{2}\mathcal{O}_{K}). $
Hence $ (K, p) $ satisfies the condition (2) in Proposition 2.2 above. So does
for (R2), too. \\
By calculating class number via PARI, for $ 5 \leq p \leq 10000 $
except $ p = 2791 \\ (h(\Q(\sqrt[3]{2791^{3}-1}))= 31876011), $ we have
$ p \nmid h(\Q(\sqrt[3]{p^{3}-1})). $ So the $ p-$Hilbert class field
$ H_{p}(K) = K \subseteq \widetilde{K}, $ i.e., $ (K, p) $ satisfies the
condition (1) in Proposition 2.2 above. Therefore by the above discussion and
Proposition 2.2, $ \Q(\sqrt[3]{p^{3}-1}) $ is $ p-$rational for all primes
$ p: 5 \leq p \leq 10000, p \neq 2791. $ On the other hand, by Gras' method
on $ p-$rationality using ray class group [G1], via calculating with
computer, one can know that $ \Q(\sqrt[3]{p^{3}-1}) $ are $ p-$rational for
all primes $ p:  p \leq 10000 $ including $ p = 2791. $  \\
Furthermore, it is easy to know that the normal closure of $ \Q(\sqrt[3]{p^{3}-1}) $
is $ K =  \Q(\sqrt[3]{p^{3}-1}, \sqrt{-3}),$ $ \text{Gal}(K / \Q) \cong S_{3}. $
So by the above discussion and Prop.3.8 of [Gr2], we can obtain that
the sextic number fields $ K =  \Q(\sqrt[3]{p^{3}-1}, \sqrt{-3}) $ are
$ p-$rational for all primes $ p : 5 \leq p \leq 10000. $ \\
More generally, we have the following conclusion: \\
Let $ K $ be a normal closure of a complex cubic number field $ k, F $ the
unique quadratic subfield of $ K,  p \geq 5 $ a prime number. Then $ K $
is $ p-$rational if and only if both $ k $ and $ F $ are $ p-$rational.
In particular, if $ k $ is a pure cubic number field, then $ F = \Q(\sqrt{-3}), $
and then $ K $ is $ p-$rational if and only if $ k $ is $ p-$rational.
\par \vskip 0.2 cm

{\bf Remark.} \ There are a series of work on the $ p-$rationality of real quadratic fields and
multi-quadratic fields, see e.g., [BM],[BR],[CLS],[Ko]. \\
It is well known that, for an odd prime $ p, $ a multi-quadratic field is $p-$rational if and only if
all of its quadratic subfields are $ p-$rational.Thus, the results on multi-quadratic fields are obtained
mainly by exploiting the $ p-$rationality of a family of quadratic fields. Some works about real quadratic
fields in the above papers are of certain types like $ \Q(\sqrt{p^2-1}),\Q(\sqrt{p^2-4}),\Q(\sqrt{p(p+ 2)}). $
These types of fields have inspired us to investigate the $ p-$rationality of the complex cubic fields
considered in this example. The common feature of them is the occurrence of a fundamental unit that can be
written down explicitly. Consequently, we can easily verify part (2) of Proposition 2.2 above.
However, for complex cubic fields, without finding a suitable explicit class number estimate, we
cannot obtain the results as complete as in the quadratic case.
\par \vskip 0.2 cm

{\bf Example 5.2.} \ $ f(x) = x^{4} - 2 x^{2} + 3 $ has four roots
$ \pm \sqrt{1 \pm \sqrt{-2}}. $ Take $ \alpha = \sqrt{1 + \sqrt{-2}}. $
Then $ k = \Q(\alpha ) $ is a pure imaginary quartic number field,
$ r_{2}(k) = 2. $ Calculating by PARI/GP, the class number $ h(k) = 1,
\varepsilon = \alpha ^{3} + \alpha ^{2} - \alpha - 2 $ is a fundamental
unit of $ k, $ and $ \{1, \alpha ,\alpha ^{2}, \alpha ^{3} \} $ is
an integral base of $ k, $ Take $ p = 5. $ It is easy to verify that
$ 5 \nmid [\mathcal{O}_{k} : \Z[\alpha]], $ and $ f(x) \text{mod}\ 5 $ is
irreducible over $ \F_{5}. $ Hence $ 5 $ is inert in $ k $ with residual
degree 4. Also by calculating, $ \varepsilon ^{5^{4} - 1} \equiv 1 + 5 \alpha
+ 15 \alpha ^{3} \not\equiv 1\ (\text{mod}\ 25 \mathcal{O}_{k}). $ Hence by
Proposition 2.2 above, $ \Q (\sqrt{1 + \sqrt{-2}}) $ is $ 5-$rational.
\par \vskip 0.2 cm

Let $ F $ be a number field, $ p $ a prime number. Let $ \widetilde{F} $
denote the composite of all $ \Z_{p}-$extension of $ F, $ and write
$ \Lambda : = \Z_{p}[[\text{Gal}(\widetilde{F}/F)]]. $ Denote by $ H(\widetilde{F}) $
the maximal unramified Abelian pro-$p$ extension of $ \widetilde{F}, $ and
write $ X(\widetilde{F}) = \text{Gal}(H(\widetilde{F})/\widetilde{F}). $ Then
$ X(\widetilde{F}) $ is a $ \Lambda -$module.
\par \vskip 0.2 cm

{\bf Cojecture 5.3}(Greenberg [Gr1]). \ {\it The $ \Lambda -$module $ X(\widetilde{F}) $
is pseudo-null. }
\par \vskip 0.2 cm

This is the well-known Greenberg's generalized conjecture (GGC, for short)(see e.g.,
[F],[It],[Mi], [T]).
\par \vskip 0.2 cm

{\bf Lemma 5.4.} \ {\it Let $ K $ be a CM field, $ K^{+} $ be its maximal real subfield,
let $ p $ be an odd prime number which splits completely in $ K. $ Let $ K_{\infty}^{+} $
be the cyclotomic $ \Z_{p}-$extension of $ K^{+}.$ Assume $ p \nmid h(K), $
where $ h(K) $ is the class number of $ K. $ Then  \\
(1) \ The Iwasawa invariants $ \lambda, \mu, \nu $ of $ K_{\infty}^{+} / K^{+} $ all vanish
and the Leopoldt's conjecture holds for $ (K^{+}, p) $ if and only if $ K^{+} $ is
$ p-$rational.  \\
(2) \ If $ K^{+} $ is $ p-$rational, then GGC holds for $ (K, p). $  }
\par \vskip 0.2 cm

{\bf Proof.} \  By definition, $ K^{+} $ is $ p-$rational if and only if
$ \text{Gal}(K_{S_{p}}^{+} / K^{+}) $ contains no $ \Z_{p}-$torsion elements and the Leopoldt's
conjecture holds for $ K^{+}. $ So by Lemma 8 of [F], conclusion (1) holds,
and then conclusion (2) follows from Theorem 1 in [F]. The proof is completed. \quad $ \Box $
\par \vskip 0.2 cm

{\bf Remark 5.5.} \ Under the assumption of Lemma 5.4 above, by the proof of Prop.2.14
of [BM], one can obtain that, $ K^{+} $ is $ p-$rational $ \Leftrightarrow $
$ K $ is $ p-$rational. So the conclusion (2) of Lemma 5.4 above can also be
stated as: If $ K $ is $ p-$rational, then GGC holds for $ (K, p). $ This partially
answers the following open question in [ABo]:   \\
Do $ p-$rational number fields satisfy Greenberg's generalized conjecture ?
\par \vskip 0.2 cm

{\bf Example 5.6.} \ Let $ p $ be a prime number, $ K = \Q(\sqrt{p^{2}-4}).$
 All such $ K $ are $ p-$rational (see [BM],Proposition 4.4). Let $ L =
K (\sqrt{-1}) = \Q(\sqrt{p^{2}-4}, \sqrt{-1}). $ We come to show that, there
are infinitely many prime number $ p $ such that GGC holds for
$ (\Q(\sqrt{p^{2}-4}, \sqrt{-1}), p). $ \\
To see this, let $ K_{1} = K, K_{2} = \Q(\sqrt{4 - p^{2}}), K_{3} =
\Q(\sqrt{-1}) $ be the three non-trivial intermediate subfields of $ L / \Q $
as the following diagram \\
We will discuss in the following steps \\
(Step1) \ If $ p \nmid h(L), $ then GGC holds for $ (L, p). $ \\
In fact, if $ p \equiv 1 (\text{mod} \ 4), $ then by our assumption, $ p $
splits completely in $ L / \Q. $ Since $ K_{1} $ is $ p-$rational and is the
maximal real subfield of $ L, $ by Lemma 5.4 above, GGC holds for $ (L, p). $
If $ p \equiv 3 (\text{mod} \ 4), $ then $ p $ splits in $ K_{2} / \Q, $
and each prime of $ K_{2} $ above $ p $ is inertia in $ L / K_{2}. $ So by
([Mi], Prop.3.B, Remark3) we know that GGC holds for $ (L, p). $ \\
(Step2) \ Now we consider the class number. By Kuroda's class number formula [Le],
$$ h(L) = \frac{1}{2} q(L) \frac{h(K_{1})h(K_{2})h(K_{3})}{h(\Q)} =
\frac{1}{2} q(L)h(K_{1})h(K_{2}),    \quad (5.6) $$
where $ q(L) = (E_{L} : E_{K_{1}}E_{K_{2}}E_{K_{3}}). $ We have $ q(L) \leq 2 $
because $ (E_{L} : W(L)E_{K}) \leq 2 $ (see [W],Thm.4.12), $ W(L) $ is the group
of roots of unity of $ L. $ On the other hand, by the relation of class number
and $ p-$rationality ([Gr2], Prop.4.1), $ p \nmid h(K_{1}). $ Therefore, by formula (5.6)
above, we have that $ p \nmid h(L) \ \Leftrightarrow \ p \nmid h(K_{2}). $ \\
 From the proof of Theorem 1 in [Ko], we know that there are infinitely many primes
$ p $ such that $ h(\Q(\sqrt{4 - p^{2}})) < p. $ In particular, such $ p $ satisfying
$ p \nmid h(\Q(\sqrt{4 - p^{2}})), $ hence GGC holds for the corresponding pair
$ (\Q(\sqrt{p^{2} - 4}, \sqrt{-1}), p). $

\par \vskip 0.2 cm

{\bf Remark 5.7.} We can prove that there are infinitely many primes
$ p $ such that GGC holds for the corresponding pair
$ (\Q(\sqrt{p^{2} - 1}, \sqrt{-1}), p)$ by the same method.
\newpage

\hspace{-0.6cm}{\bf 6.\ Tables }

      \par \vskip 0.6 cm

\begin{table}[!ht]
	\center
\caption{Complex cubic number field}
	\begin{tabular}{|c|c|c|c|}\hline
		\multicolumn{1}{|c|}{\multirow{2}{*}{$ T(|d_{K}|<3000) $}} & \multicolumn{3}{c|}{$5\leq p\leq 100$} \\ \cline{2-4}
		&$p|h(K)$ &$\varepsilon^{p^f-1}\equiv 1\ (\ mod\ \mathfrak{p}^{e+1}),$ $\forall\mathfrak{p}|p$ &not $p$-rational \\ \hline
		$x^3-x^2+x-9$&$\emptyset$&13 & 13\\ \hline
	$ x^3-x^2+5*x+1$&$\emptyset$ &17&17\\ \hline
	$x^3-x^2-2*x+6 $&$\emptyset$ &5&5\\ \hline
	$ x^3-x^2+x+5$&$\emptyset$&5&5 \\ \hline
	$x^3-x^2+5*x+2$ &$\emptyset$&11&11\\ \hline
	$x^3-6*x-12$&$\emptyset$&5&5\\ \hline
	$x^3-x^2-x+13$&$\emptyset$&5&5\\ \hline
	$x^3-x^2-x-6$ &$\emptyset$&11&11\\ \hline
	$x^3-x^2+5*x+11$&$\emptyset$&11&11\\ \hline
	$x^3-x^2+7*x-2$&$\emptyset$&5&5\\ \hline
	$x^3-8*x-11$&$\emptyset$&5&5\\ \hline
	$x^3-x^2-4*x+9$&$\emptyset$&19&19\\ \hline
	$x^3-x^2+7*x-6$&5&$\emptyset$ &5?\\ \hline
	$x^3-x^2+x+15$&5& $\emptyset$&5?\\ \hline
	$x^3-x^2+x-24$&$\emptyset$&5&5\\ \hline
	$x^3-x^2+4*x-9$&$\emptyset$&13&13\\ \hline
	$x^3-x^2-6*x-16$&$\emptyset$&5&5\\ \hline
	$x^3+10*x-12$&5&$\emptyset$&5?\\ \hline
	$x^3-x^2+10*x-16$&5&$\emptyset$&5?\\ \hline
	$x^3-26$&$\emptyset$&11&11\\ \hline
	$x^3-x^2-8*x-10$&$\emptyset$&5&5\\ \hline
	$x^3-x^2-x-26$&$\emptyset$&61&61\\ \hline
	$x^3-x^2+13*x-1$&5&$\emptyset$&5?\\ \hline
	$x^3-x^2-3*x-17$&$\emptyset$&5&5\\ \hline
	$x^3-x^2+7*x-19$&$\emptyset$&31&31\\ \hline
	$x^3-x^2-11*x+21$&$\emptyset$&11&11\\ \hline
	$x^3-11*x-17$&$\emptyset$&5&5\\ \hline
	$x^3-x^2+6*x-10$&$\emptyset$&23&23\\ \hline
	$x^3-x^2-10*x-20$&5&$\emptyset$&5? \\ \hline
   $x^3-x^2-11*x-21$&$\emptyset$&7&7\\ \hline
   $x^3-2*x-20$&5&$\emptyset$&5?\\ \hline
   $x^3+2*x-10$&$\emptyset$&7&7\\ \hline
   $x^3+4*x-20$&$\emptyset$&13&13\\ \hline
   $x^3-x^2+5*x-32$&5&$\emptyset$&5?\\ \hline
   $x^3-x^2+9*x-21$&7&$\emptyset$&7?\\ \hline	

\end{tabular}
\end{table}

\begin{table}[!ht]
\center
\caption{Pure imaginary quartic number field}
\begin{tabular}{|c|c|c|c|c|}\hline
		\multicolumn{1}{|c|}{\multirow{2}{*}{$ T (d_{K}<3000) $}} & \multicolumn{4}{c|}{$5\leq p\leq 100$} \\ \cline{2-5}
		&5 totally ramified &$p|h(K)$ &$\varepsilon^{p^f-1}\equiv 1\ (\ mod\ \mathfrak{p}^{e+1}),$ $\forall\mathfrak{p}|p$ &not $p$-rational \\ \hline
		$x^4-x^3+x^2-x+1$&Yes&$\emptyset$&$\emptyset$&-\\ \hline
		$x^4+1$& No&$\emptyset$&13,31 & 13,31\\ \hline
		$x^4-2*x^2+4$&No &$\emptyset$&7&7\\ \hline
		$x^4+2*x^2+4$& No&$\emptyset$&13,31&13,31\\ \hline
		$x^4-2*x^3-2*x+5$&No &$\emptyset$&11&11\\ \hline
		$x^4-x^3-4*x^2+4*x+7$&No &$\emptyset$&23&23\\ \hline
		$x^4-2*x^3+5*x^2-4*x+2$&No&$\emptyset$&13,31&13,31\\ \hline
		$x^4-x^3-2*x^2-3*x+9$&No&$\emptyset$&29,37&29,37\\ \hline
		$x^4-2*x^3-4*x^2+5*x+7$&No&$\emptyset$&5&5\\ \hline
		$x^4-2*x^3-3*x^2+4*x+5$&No&$\emptyset$&11&11\\ \hline	
		 $x^4+4*x^2+2$&No&$\emptyset$&13,31&13,31\\ \hline
		$x^4+9$&No&$\emptyset$&7&7\\ \hline
		
	\end{tabular}
\end{table}

\newpage

In the above two tables, we use LMFDB data: In Table 1, there are totally $ 419 $
number of complex cubic number fields $ K $ of which the absolute values of discriminants
$ < 3000. $ We take all the $ 23 $ number of prime numbers $ p: 5 \leq p \leq 100. $
For the $ 419 \times 23 = 9637 $ pair of $ (K, p), $ by our criteria Proposition 2.2,
except the cases of the above table, all the remainder are $ p-$rational. In the table, we
use a polynomial $ T $ of degree $ 3 $ to indicate the cubic number field which
is generated by a root of $ T. $ For the case that $ p \mid h(K) $ with
$\varepsilon^{p^f-1}\equiv 1\ (\ mod\ \mathfrak{p}^{e+1}),\ \forall\mathfrak{p}|p,$
to determine the $ p-$rationality, one need to further determine whether
$ H_{p} (K) \subset \widetilde{K}. $ At present, we don't know whether there are
such effective algorithms. \\
In Table 2, there are totally $ 147 $ number of pure imaginary quartic number fields
$ K $ of which the absolute values of discriminants $ < 3000. $ We take all the $ 23 $
number of prime numbers $ p: 5 \leq p \leq 100. $
For the $ 147 \times 23 = 3381 $ pair of $ (K, p), $ by our criteria Proposition 2.2,
except the cases of the above table, all the remainder are $ p-$rational. In the table, we
use a polynomial $ T $ of degree $ 4 $ to indicate the quartic number field which
is generated by a root of $ T. $ We need to point out, at present, for the case that
$ 5 $ is totally ramified in the quartic number field, we can not use the result of
Proposition 2.2 above to determine the $ 5-$rationality. \\
For the determined cases, the calculated results in our tables
are consistent with the ones done by the algorithm of Gras in [G1].
\par  \vskip 0.3cm

{ \bf Acknowledgments } \ We would like to thank Georges Gras
for sending us his several papers and giving helpful suggestions. We would like to thank
Zouhair Boughadi for informing us the relation of an open question in [ABo] and
our Lemma 5.4 as in Remark 5.5. We would like to thank Ruiyuan He and Qinhao Li for their
careful reading of this paper and helpful suggestions. We would like to thank the referee
for helpful suggestions and comments.

\par  \vskip 1 cm
\hspace{-0.8cm} {\bf References }
\begin{description}

\item[[AB]]  J. Assim and Z. Bouazzaoui, Half-integral weight modular forms and real
quadratic $ p$-rational fields. Funct. Approx. Comment. Math., 63(2):201-213,2020.
\item[[ABo]]  J. Assim and Z. Boughadi, On Greenberg's generalized conjecture, J. Number Theory,
242: 576-602, 2023.
\item[[BM]]  Y.Benmerieme, A.Movahhedi, Multi-quadratic $ p$-rational number fields. J. Pure Appl.
Algebra,225(9): 17, 2021.
\item[[BR]] R.Barbulescu, J.Ray, Numerical verification of the Cohen-Lenstra-Martinet heuristics and
of Greenberg's $p-$rationality conjecture, J.Th$ \acute{e}$or.Nombres Bordeaux 32(2020), no.1, 159-177.
\item[[CLS]]  J. Chattopadhyay, H. Laxmi, A. Saikia, On the $ p-$rationality of consecutive
quadratic fields, J. Number Theory, 248: 14-26, 2023.
\item[[F]]  S. Fujii, On Greenberg's generalized conjecture for CM-fields. J. Reine Angew. Math.,
731:259-278,2017.
\item[[G1]] G. Gras, On $ p-$rationality of number fields. Applications---PARI/GP programs.
In Publications math\'{e}matiques de {B}esan\c{c}on. {A}lg\`ebre et th\'{e}orie des nombres.
2019/2, volume 2019/2 of Publ. Math. Besan\c{c}on Alg\`ebre Th\'{e}orie Nr., 29-51. Presses
Univ. Franche-Comt\'{e}, Besan\c{c}on, 2019.
\item[[G2]] G. Gras, les $ \theta-$r$\acute{e}$gulateurs locaux d'un nombre alg$\acute{e}$braique:
conjectures $ p-$adiques, Can. J. Math. 68(3), 2016: 571-624.
\item[[Gr1]]  R. Greenberg, Iwasawa theory--past and present. In Class field theory--its
centenary and prospect (Tokyo, 1998), volume 30 of Adv. Stud. Pure Math., pages 335-385.
Math.Soc.Japan, Tokyo,2001.
\item[[Gr2]]  R. Greenberg, Galois representations with open image. Ann. Math. Qu\'{e}.,
40(1):83-119,2016.
\item[[Ic]]  H. Ichimura, A note on Greenberg's conjecture and the $ abc $ conjecture,
Proc. Amer. Math. Soc. 126(1998), no.5, 1315-1320.
\item[[It]]  T. Itoh, On multiple $ \Bbb Z_p$-extensions of imaginary abelian
quartic fields. J. Number Theory, 131(1):59-66,2011.
\item[[Ka]]  T. Kataoka, On Greenberg's generalized conjecture for complex cubic fields,
Int. J. Number Theory, 133:619-631,2017.
\item[[Ko]]  J. Koperecz, Triquadratic $ p$-rational fields. J. Number Theory, 242:402-408,2023.
\item[[Le]]  F. Lemmermeyer, Kuroda's class number formula. Acta Arith., 66(3):245-260,1994.
\item[[Li]]  D. Lim, On $ p$-rationality of $\Bbb Q(\zeta_{2l+1})^+ $ for
Sophie Germain primes $ l. $ J. Number Theory, 231:378-400,2022.
\item[[Mi]]  J. V. Minardi, Iwasawa modules for $ \Bbb Z^d_p $-extensions of algebraic
number fields. ProQuest LLC, Ann Arbor, MI,1986. Thesis(Ph.D.)-University of Washington.
\item[[Mov1]]  A. Movahhedi, Sur les $ p$-extensions des corps $ p$-rationnels.,
PhD. Thesis, 1988.
\item[[Mov2]]  A. Movahhedi, Sur les $ p$-extensions des corps $ p$-rationnels. Math. Nachr.,
149:163-176,1990.
\item[[MR]] C. Maire, M. Rougnant, A note on $ p-$rational fields and the $ abc-$conjecture.
Proc. Amer. Math. Soc. 148(2020), no.8, 3263-3271.
\item[[N]] J. Neukirch, Algebraic number theory, volume 322 of Grundlehren der
mathematischen Wissenschaften [Fundamental Principles of Mathematical Sciences].
Springer-Verlag, Berlin,1999. Translated from the 1992 German original and with a note by
Norbert Schappacher, With a foreword by G. Harder.
\item[[NG]] T. Nguyen Quang Do, On Greenberg's generalized conjecture for families of number
fields. arXiv: 2505.07529v1.
\item[[S]] J.-P. Serre, Local fields, volume 67 of Graduate Texts in Mathematics.
Springer-Verlag, New York-Berlin,1979.Translated from the French by Marvin Jay Greenberg.
\item[[Sil]] J. H. Silverman, Wieferich's criterion and the $ abc-$conjecture.
J. Number Theory, 30 (1988), no. 2, 226-237.
\item[[St]] H.-J. Stender, \"Uber die Grundeinheit f\"{u}r spezielle unendliche Klassen
reiner kubischer Zahlk\"{o}rper. Abh. Math. Sem. Univ. Hamburg, 33:203-215,1969.
\item[[T]] N. Takahashi, On Greenberg's generalized conjecture for imaginary quartic fields,
Int. J. Number Theory, 17 (05) (2021), 1163-1173.
\item[[W]] L. C. Washington, Introduction to cyclotomic fields, volume 83 of Graduate
Texts in Mathematics. Springer-Verlag, New York,second edition,1997.

\end{description}

\end{document}